# A Comparison of the Carlitz and Digit Derivatives Bases in Function Field Arithmetic

Sangtae Jeong

Department of Mathematics, Korea University, Seoul, Korea 136-701

stj@semi.korea.ac.kr

**Abstract**: We compare several properties and constructions of the Carlitz polynomials and digit derivatives for continuous functions on $\mathbf{F}_q[[T]]$. In particular, we show a close relation between them as orthonormal bases. Moreover, parallel to Carlitz's coefficient formula, we give the closed formula for the expansion coefficients in terms of the digit derivatives.

## 1 Introduction

As an analogue of Mahler expansion for function fields, Wagner [W1] and Goss [Go1] proved that Carlitz polynomials $\{G_j\}_{j\geq 0}$ are an orthonormal basis of the Banach space of continuous functions from $O := \mathbf{F}_q[[T]]$ to $K := \mathbf{F}_q((T))$, denoted $C(O, K)$, which is topologized by the sup-norm. Unlike the case of $p$-adic integers, the space $C(O, K)$ has a subspace of continuous $\mathbf{F}_q$-linear functions from $O$ to $K$, denoted $LC(O, K)$. It is also shown by them [W2], [Go1] that $LC(O, K)$ has an orthonormal basis of Carlitz $\mathbf{F}_q$-linear polynomials, $\{E_n(x)\}_{n\geq 0}$. As it is shown by Carlitz [C1], the Carlitz polynomials $\{G_j\}_{j\geq 0}$ are formed out of the orthonormal basis of the subspace, $\{E_n(x)\}_{n\geq 0}$ via the $q$-adic expansion of $j$.

Recently, B. Snyder [Sn] and the author [J] showed with different arguments that the certain $\mathbf{F}_q$-linear hyper-differential operators $\{\mathcal{D}_n\}_{n\geq 0}$ are an orthonormal basis of $LC(O, K)$. The digit derivatives $\{\mathsf{D}_j\}_{j\geq 0}$ are formed from hyper-differential operators $\{\mathcal{D}_n\}_{n\geq 0}$ via the $q$-adic expansion of $j$ in the same way that Carlitz did with the $E_n(x)$. Then we see by the digit principle(due to K. Conrad [C2]) that the digit derivatives $\{\mathsf{D}_j\}_{j\geq 0}$ as well as Carlitz polynomials $\{G_j\}_{j\geq 0}$ are orthonormal bases of $C(O, K)$. Indeed, any orthonormal basis of $LC(O, K)$ can extend to an orthonormal basis of $C(O, K)$ via the $q$-adic expansion of non-negative integers. His method verifies a sequence of $q$-adic extensions is an orthonormal basis of the whole space but gives no information about the expansion coefficients.

In this note, from the concrete point of view, the focus will be on explicit expansions of functions $f \in C(O, K)$ in terms of the Carlitz polynomials and digit derivatives as well as on those for functions $f \in LC(O, K)$ in terms of Carlitz linear polynomials and hyper-differential operators. We explicitly compare the expansion coefficients for any function $f \in LC(O, K)$ in terms of two canonical orthonormal bases for the subspace. We then show the two bases are essentially equivalent, which means here that Carlitz linear polynomials being an orthonormal basis for $LC(O, K)$ can be deduced from knowing that hyper-differential operators are an orthonormal basis and vice versa. The equivalence relation between the two bases for the subspace $LC(O, K)$ gives us that both Carlitz polynomials $\{G_j\}_{j\geq 0}$ and digit derivatives $\{\mathsf{D}_j\}_{j\geq 0}$ have the same images in the reduced space $C(O, \mathbf{F}_q)$. Thus we see with no further work that the $q$-adic extensions are orthonormal bases of $C(O, K)$.

As for the coefficients for expansion of functions in $C(O, K)$, Carlitz [C2] could recover them in terms of the Carlitz polynomials in the closed formula. Parallel to this formula, we will also obtain the closed formula for the expansion coefficients of functions in $C(O, K)$ in terms of the digit derivatives by following Yang's arguments [Y]. Indeed this is based on the close similarity between the orthogonality properties for two canonical orthonormal bases of $C(O, K)$.

Finally we show for any non-negative integer $m$ that the two sets of $q^m$ powers of hyper-differential operators and of Carlitz $\mathbf{F}_q$-linear polynomials are also orthonormal bases of $LC(O, K)$ and then obtain coefficient formulas of expansion of any continuous linear functions with respect to $q^m$ powers of hyper-differential operators in the closed form related to the Carlitz difference operators.

**Acknowledgements.** The author wishes to thank K. Conrad and D. Goss for their invaluable comments on the earlier draft of the paper.

## 2 Two canonical orthonormal bases

Let $K := \mathbf{F}_q((T))$ be the field of formal power series over a finite field $\mathbf{F}_q$ of $q$ elements. Elements of $K$ are written uniquely as a Laurent series expansion with coefficients in $\mathbf{F}_q$. That is, for non-zero $x \in K$, $x = \sum_{i=-\infty}^{\infty} \alpha_i T^i$, where $\alpha_i = 0$ for all but finite many $i \leq 0$. Then the valuation $v$ of $K$ is given by $v(x)$ = the smallest integer $i$ for which $\alpha_i \neq 0$ in the expansion of $x$. The absolute value $|x|$ on $K$ is defined as $|x| = q^{-v(x)}$ for $x \in K^*$. By convention $|0| = 0$. Then $K$ is a complete non-archimedean field of a non-trivial absolute value $|\cdot|$.

Let $O := \mathbf{F}_q[[T]]$ be the ring of integers in $K$ consisting of formal power series of the form $\sum_{i=0}^{\infty} \alpha_i T^i$. Recall that $C(O, K)$ denotes the $K$-Banach space of all continuous functions $f : O \to K$, equipped with the sup-norm $\|f\| = \max_{t \in O}\{|f(t)|\}$ and $LC(O, K)$ is the subspace of continuous $\mathbf{F}_q$-linear functions from $O$ to $K$.

Let $e_0(x) = x$ and for $n \geq 1$ let
$$e_n(x) = \prod_{\substack{m \in \mathbf{F}_q[T] \\ \deg(m) < n}} (x - m),$$

Let $F_0 = L_0 = 1$ and for $n \geq 1$ let
$$F_n = [n][n-1]^q \cdots [1]^{q^{n-1}}, \quad L_n = [n][n-1] \cdots [1],$$

where $[n] = T^{q^n} - T$. It is well-known [C1],[Go2] that $e_n(x)$ is an $\mathbf{F}_q$-linear polynomial of degree $q^n$ with coefficients in $\mathbf{F}_q[T]$, as it has an expansion:
$$e_n(x) = \sum_{i=0}^{n} (-1)^{n-i} \frac{F_n}{F_i L_{n-i}^{q^i}} x^{q^i}.$$

Another set of functions in $LC(O, K)$ is the hyper-differential operators, namely Hasse derivations $\mathcal{D}_n, n \geq 0$, defined by
$$\mathcal{D}_n(\sum a_i T^i) = \sum \binom{i}{n} a_i T^{i-n}.$$

It is known in [V] that $\mathcal{D}_n$ are continuous $\mathbf{F}_q$-linear functions on $K$ (hence on $O$). For later use, recall here the properties of hyper-differential operators $\mathcal{D}_n$;

(1) $\mathcal{D}_0 = id$

(2) $\mathcal{D}_n(xy) = \sum_{i=0}^{n} \mathcal{D}_i(x)\mathcal{D}_{n-i}(y)$

(3) $\mathcal{D}_n(\mathcal{D}_m(x)) = \binom{n+m}{m} \mathcal{D}_{n+m}(x) \quad \text{for all} \quad x, y \in K$.

The collection of additive operators satisfying the three properties above is called iterative derivations and for a detailed study on them, see [Sch].

**Definition 1** (1) *Put $E_n(x) = e_n(x)/F_n$ for any integer $n > 0$ and $E_0(x) = x$.*

(2) *For the $q$-adic expansion of $j \geq 0$, which is given by*
$$j = \alpha_0 + \alpha_1 q + \cdots + \alpha_s q^s$$

*with $0 \leq \alpha_i < q$, Put*
$$G_j(x) := \prod_{n=0}^{s} E_n^{\alpha_n}(x), j \geq 1; G_0(x) = 1,$$

and
$$G'_j(x) := \prod_{n=0}^{s} G'_{\alpha_n q^n}(x),$$

where

$$G'_{\alpha q^n}(x) = \begin{cases} E_n^\alpha(x) & if 0 \leq \alpha < q-1 \\ E_n^\alpha(x) - 1 & if \alpha = q-1. \end{cases}$$

(3) *Put*

$$\mathsf{D}_j(x) := \prod_{n=0}^{s} \mathcal{D}_n^{\alpha_n}(x), j \geq 1; \mathsf{D}_0(x) = 1,$$

and

$$\mathsf{D}'_j(x) := \prod_{n=0}^{s} \mathsf{D}'_{\alpha_n q^n}(x),$$

where

$$\mathsf{D}'_{\alpha q^n}(x) = \begin{cases} \mathcal{D}_n^\alpha(x) & if 0 \leq \alpha < q-1 \\ \mathcal{D}_n^\alpha(x) - 1 & if \alpha = q-1. \end{cases}$$

Thus in particular, we see that

$$G_{\alpha q^n}(x) = E_n^\alpha(x), \mathsf{D}_{\alpha q^n}(x) = \mathcal{D}_n^\alpha(x) \quad \text{for} \ \ 0 \leq \alpha < q.$$

In what follows, we will call $G_j, \mathsf{D}_j$ Carlitz polynomials, digit derivatives, respectively. Observe the difference between the polynomials $G_j(x)$ and $G'_j(x)$ defined here and Carlitz's original polynomials $G_j(x)/g_j$ and $G'_j(x)/g_j$ defined in [C2], where $g_j := \prod_{n=0}^{s} F_n^{\alpha_n}$ is an analogue of the classical factorial. These $G_j(x)$ and $G'_j(x)$ are polynomials of degree $j$ and satisfy various identities such as the orthogonality property, which will be discussed in section 4. Moreover it is known in [C2] that $G_j(x)$ and $G'_j(x)$ are integral-valued polynomials , which means that $G_j(m)$ and $G'_j(m)$ belong to $\mathbf{F}_q[T]$ for any $m \in \mathbf{F}_q[T]$. Parallel to this it is easy to see that $\mathsf{D}_j(x)$ and $\mathsf{D}'_j(x)$ are also integral-valued functions.

**Definition 2** *Let $K$ be a local field, and $E$ be a $K$-Banach space equipped with the usual sup-norm. We say that a sequence $\{f_n\}_{n \geq 0}$ in $E$ is an orthonormal basis for $E$ if and only if the following two conditions are satisfied:*

(1) *every $f \in E$ can be expanded uniquely as $f = \sum_{n=0}^{\infty} a_n f_n$, with $a_n \in K \to 0$.*
(2) *The sup-norm of $f$ is given by $\|f\| = max\{|a_n|\}$.*

Wagner [W2] first established an analogue of Mahler's expansion for continuous $\mathbf{F}_q$-linear operators on $O$ in terms of Carlitz $\mathbf{F}_q$-linear polynomials, $E_n(x)$.

**Theorem 1** *The set of Carlitz $\mathbf{F}_q$-linear polynomials, $\{E_n(x)\}_{n \geq 0}$ is an orthonormal basis for $LC(O, K)$. That is, every continuous $\mathbf{F}_q$-linear function $f$ can be expanded uniquely as $f = \sum_{n=0}^{\infty} a_n E_n(x)$ with $a_n \to 0$, and the sup-norm of $f$ is given by $\|f\| = max\{|a_n|\}$. Moreover the coefficients can be recovered by the formula:*

$$a_n = (\Delta^{(n)} f)(1),$$

*where the Carlitz difference operators $\Delta^{(n)}, n \geq 0$ are defined recursively by*

$$(\Delta^{(n)} f)(x) = \Delta^{(n-1)} f(Tx) - T^{q^{n-1}} \Delta^{(n-1)} f(x), n \geq 1; \Delta^0 = id.$$

Voloch [V] established the following explicit expansion of $\mathcal{D}_m$ in terms of Carlitz linear polynomials $E_n(x)$.

**Proposition 1** *For any $x \in O, \mathcal{D}_m(x) = \sum_{n \geq m}^{\infty} A_{n,m} E_n$, where $A_{n,1} = (-1)^{n-1} L_{n-1}$ and for $m > 1$,*

$$A_{n,m} = (-1)^{n+m} L_{n-1} \sum_{0 < i_1 \cdots < i_{m-1} < n} \frac{1}{[i_1][i_2]\cdots[i_{m-1}]}.$$

Note that $A_{n,n} = 1$ for each non-negative integer $n$.

As a result analogous to Wagner's, Snyder [Sn] and the present author [J] independently established the representation of continuous $\mathbf{F}_q$-linear functions with respect to hyper-differential operators in the following.

**Theorem 2** *The set of hyper-differential operators, $\{\mathcal{D}_n\}_{n \geq 0}$ is an orthonormal basis for $LC(O, K)$. That is, every $f$ can be expanded uniquely as $f = \sum_{n=0}^{\infty} b_n \mathcal{D}_n$ with $b_n \to 0$, and the sup-norm of $f$ is given by $\|f\| = \max\{|b_n|\}$. Moreover the coefficients can be recovered by iterating the Carlitz difference operator $\Delta$:*

$$b_n = (\Delta^n f)(1) = \sum_{i=0}^{n} (-1)^{n-i} f(T^i) \mathcal{D}_i(T^n),$$

where $\Delta := \Delta^{(1)}$ is defined above.

By taking $E_n(x)$ for $f$ in Theorem 2 the following is immediate.

**Corollary 1** *For $x \in O$, $E_n(x) = \sum_{m=0}^{\infty} B_{m,n} \mathcal{D}_m(x)$, where*

$$B_{m,n} = \sum_{i=0}^{m} (-1)^{m-i} \mathcal{D}_i(T^m) E_n(T^i).$$

We also note that $B_{m,n} = 0$ for $m < n$, and $B_{n,n} = 1$ as $E_n(T^n) = 1$. Notice that the expansion in Corollary 1 is the inversion formula to the Voloch's expansion stated in Proposition 1.

In order to establish Theorem 2, Snyder follows the arguments of Wagner but the proof of the author is based on Theorem 50.7 [S] on orthonormal bases of Banach spaces over local fields. In fact, both proofs depend on Wagner's result (Theorem 1) and use Voloch's result. On the other hand, Conrad [Co2] proved the same result, independent of Wagner by simply using a criterion on orthonormal bases for Banach spaces over local fields and did not find the coefficient formulas. But the formula for the coefficients is immediate, as $\Delta$ acts on $\mathcal{D}_n$ by shifting, namely $\Delta \mathcal{D}_n = \mathcal{D}_{n-1}$. More importantly, there he also established the digit principle, which will be used in the proof of Theorem 6. This principle is the reason we are interested in bases for subspaces of $C(O, K)$. Before we formulate it, recall the $q$-adic expansion of an arbitrary integer $j \geq 0$ with

$$j = \alpha_0 + \alpha_1 q + \cdots + \alpha_s q^s,$$

where $0 \leq \alpha_i < q$. For a given orthonormal basis $\{e_n\}_{n \geq 0}$ of $LC(O, K)$, we associate the function

$$E_j := e_0^{\alpha_0} e_1^{\alpha_2} \cdots e_s^{\alpha_s}.$$

It turns out that this collection of functions, $E_j$ for all $j$, is an orthonormal basis of $C(O, K)$.

**Theorem 3 (Digit Principle)**
*If $\{e_n\}_{n \geq 0}$ is an orthonormal basis of $LC(O, K)$, then $\{E_j\}_{j \geq 0}$ is an orthonormal basis of $C(O, K)$.*

**Proof.** See [Co2] for proof. □

By applying Conrad's digit principle to Carlitz linear polynomials and hyper-differential operators the following theorem is immediate. On the other hand, his principle only verifies a sequence of $q$-adic extensions is an orthonormal basis but gives no information about the coefficients. But using interpolations of Carlitz (linear) polynomials, Carlitz [C2],[W2] could recover the expansion coefficients of $f \in C(O, K)$ in terms of Carlitz polynomials in the closed form. Then Yang [Y] also proved the same closed formula of the coefficients in an elementary way by using the orthogonality property of $G_j(x)$ and $G'_j(x)$. To sum it up, we have

**Theorem 4** (1) *Both* $\{G_j(x)\}_{j \geq 0}$ *and* $\{\mathsf{D}_j(x)\}_{j \geq 0}$ *are orthonormal bases of* $C(O, K)$.

(2) *Write* $f \in C(O, K)$ *as* $f = \sum_{j \geq 0} A_j G_j(x)$. *Then* $A_j$ *can be recovered by the formula: For any integer $n$ such that $q^n > j$,*

$$A_j = (-1)^n \sum_{\substack{m \in \mathbf{F}_q[T] \\ deg(m) < n}} G'_{q^n - 1 - j}(m) f(m).$$

The formula for coefficients $A_j$ gives the following.

**Corollary 2** *Let* $f(x) = \sum_{j \geq 0} A_j G_j(x)$ *be a continuous function from $O$ to $K$. Then $f(x) \in C(O, O)$ if and only if $\{A_j\}_{j \geq 0} \subset O$.*

We close this section with another proof of Theorem 2 by using a useful criterion([Se] Lemme I) as Conrad did with hyper-differential operators. This is really easy to see and does not depend on Wagner's arguments. A necessary and sufficient condition for Carlitz $\mathbf{F}_q$-linear polynomials, $\{E_n(x)\}_{n \geq 0}$, to be an orthonormal basis for $LC(O, K)$ is that

(a) $E_n(x)$ maps $O$ into itself.

(b) their reduced functions, $E_n(x) \mod (T)$, denoted $\overline{E_n}$, are a basis of $LC(O, \mathbf{F}_q)$ in an algebraic sense as an $\mathbf{F}_q$-vector space.

For Part (a) it follows from Carlitz's result that $E_n(x)$ is an integral-valued $\mathbf{F}_q$-linear operator [C2](Lemma 1). To check Part (b) it is enough to show that for any $n \geq 1, \overline{E}_0, \overline{E}_1 \cdots \overline{E}_{n-1}$ is a basis of the linear space of $\mathbf{F}_q[T]/T^n$ into $\mathbf{F}_q$. By reason of cardinality it suffices to show that they are linearly independent over $\mathbf{F}_q$. Suppose that

$$\alpha_0 \overline{E}_0(x) + \alpha_1 \overline{E}_1(x) + \cdots + \alpha_{n-1} \overline{E}_{n-1}(x) = 0.$$

Plugging $1, T, \cdots T^{n-1}$ into $x$ successively gives us $\alpha_i = 0$ for $0 \leq i \leq n-1$, as $E_i(T^i) = 1$ for all $i$.

## 3 Relation between two bases

We now know that $LC(O, K)$ has two canonical sets of orthonormal bases, Carlitz $\mathbf{F}_q$-linear polynomials and hyper-differential operators. In this section, we would like to show by a simple and different method that Carlitz $\mathbf{F}_q$-linear polynomials being an orthonormal basis for $LC(O, K)$ can be derived from knowing that hyper-differential operators are an orthonormal basis and vice versa. Then we will see by a direct calculation that both Carlitz functions $\{G_j\}_{j \geq 0}$ and digit derivatives $\{\mathsf{D}_j\}_{j \geq 0}$ have the same images in the reduced space $C(O, \mathbf{F}_q)$. To this end, we first introduce a useful lemma as an easy criterion of orthonormal basis for Banach spaces over complete non-archimedean fields, which will be used in many places throughout the paper.

**Lemma 1** *Let $K$ be a complete non-archimedean field with a nontrivial absolute value and $E$ be a $K$-Banach space with an orthonormal basis $\{e_n\}_{n \geq 0}$. If $f_n \in E$ with $\sup_{n \geq 0} \|e_n - f_n\| < 1$ then $\{f_n\}_{n \geq 0}$ is an orthonormal basis of $E$.*

**Proof.** It follows immediately from Serre's criterion ([Se] Lemme I) as we easily deduce from hypothesis that $\|f_n\| = 1$ and the reduced functions for $f_n$ and $e_n$ are identically same . Alternatively we refer the reader to [Co1] Lemma 3.2 for a concrete proof. □

**Theorem 5** $\{E_n(x)\}_{n\geq 0}$ *is an orthonormal basis of $LC(O,K)$ if and only if $\{\mathcal{D}_n\}_{n\geq 0}$ is an orthonormal basis of $LC(O,K)$.*

**Proof.** We will give two proofs.
Assume that $\{E_n(x)\}_{n\geq 0}$ is an orthonormal basis of $LC(O,K)$. Then from the expansion for $\mathcal{D}_m$ in Proposition 1 and remark to $A_{m,m} = 1$, we calculate for each $m \geq 0$,

$$\|\mathcal{D}_m - E_m\| = \|\sum_{n=m+1}^{\infty} A_{n,m} E_n(x)\| = \max_{n \geq m+1}\{|A_{n,m}|\} \leq 1/q < 1,$$

as $v(A_{n,m}) \geq n - m$ for $n \geq m + 1$. Now the result follows from Lemma 1 above.

Conversely, assuming that $\{\mathcal{D}_n\}_{n\geq 0}$ is an orthonormal basis of $LC(O,K)$, we then use the expansion for $E_n(x)$ in Corollary 1. By recalling remarks to $B_{m,n}$, we get for each $n \geq 0$,

$$\|E_n - \mathcal{D}_n\| = \|\sum_{m=n+1}^{\infty} B_{m,n} \mathcal{D}_m\| = \max_{m \geq n+1}\{|B_{m,n}|\}.$$

Since $E_n(x)$ is an integral-valued polynomial [C2](Lemma 1) and $T$ divides $\mathcal{D}_i(T^m)$ for $i < m$, we need to show that $T$ divides $E_n(T^m)$ for $m \geq n+1$ in the formula for $B_{m,n}$ given in Corollary 1. It follows by repeatedly applying the Carlitz's identity [C2](8.2);

$$\text{For any integer } m, E_n(T^{m+1}) = TE_n(T^m) + E_{n-1}^q(T^m).$$

Hence for $m \geq n+1, |B_{m,n}| \leq 1/q < 1$. So $\|E_n - \mathcal{D}_n\| \leq 1/q < 1$. Thus the result again comes from Lemma 1.

We now give an alternate proof, which is suggested by Conrad. It involves more directly the original definitions of the functions $\mathcal{D}_n(x)$ and $E_n(x)$, and it does not resort to a link between $\mathcal{D}_n$ and $E_n$, as it is shown above. We will show for all $n$ that the reduced functions for $\mathcal{D}_n$ and $E_n$, as maps from $O$ to $\mathbf{F}_q$, coincide. Therefore we are done by Lemma 1.

Both $\overline{E_n}$ and $\overline{\mathcal{D}_n}$ vanish on the ideal $T^{n+1}$, so they define functions from $O/T^{n+1}$ to $\mathbf{F}_q$. Each vanishes at $T^i$ for $i < n$ and equals 1 for $T^n$. So by linearity they are equal. □

Then from Theorem 5 invoking either Theorem 1 or 2 gives the immediate corollary.

**Corollary 3** *Both $\{E_n(x)\}_{n\geq 0}$ and $\{\mathcal{D}_n\}_{n\geq 0}$ are orthonormal bases of $LC(O,K)$.*

**Theorem 6** *Carlitz functions $\{G_j\}_{j\geq 0}$ are an orthonormal basis of $C(O,K)$ if and only if digit derivatives $\{\mathsf{D}_j\}_{j\geq 0}$ are an orthonormal basis of $C(O,K)$.*

**Proof.** It follows from Theorem 5 and the digit principle. □

**Theorem 7** *(1) Both $\{G_j(x)\}_{j\geq 0}$ and $\{\mathsf{D}_j(x)\}_{j\geq 0}$ are orthonormal bases of $C(O,K)$.*
*(2) Write $f \in C(O,K)$ as $f = \sum_{j\geq 0} B_j \mathsf{D}_j(x)$. Then $B_j$ can be recovered by the formula: For any integer $n$ such that $q^n > j$,*

$$B_j = (-1)^n \sum_{\substack{m \in \mathbf{F}_q[T] \\ deg(m) < n}} \mathsf{D}'_{q^n-1-j}(m) f(m).$$

**Proof.** The first assertion follows from Theorem 6 by referring to either [W1] or [Go1]. Or it follows by applying the digit principle to Corollary 3. The proof of the second assertion begins by invoking the orthogonality property (Theorem 9(1) below) of $\mathsf{D}_j(x)$ and $\mathsf{D}'_j(x)$. For any integer $n$ such that $q^n > j$, we have

$$(-1)^n \sum_{deg(m)<n} \mathsf{D}'_{q^n-1-j}(m) f(m) = \sum_{j=0}^{\infty} B_j (-1)^n \sum_{deg(m)<n} \mathsf{D}'_{q^n-1-n}(m) \mathsf{D}_j(m).$$

The orthogonality property(Theorem 9(1)) yields the desired formula for $B_j$. □

The formula for coefficients $B_j$ yields the following.

**Corollary 4** *Let $f(x) = \sum_{j\geq 0} B_j \mathsf{D}_j(x)$ be a continuous function from $O$ to $K$. Then $f(x) \in C(O,O)$ if and only if $\{B_j\}_{j\geq 0} \subset O$.*

Recall that $K = \mathbf{F}_q((T))$. Let $T^{1/q^m}$ be a $q^m$th root of $T$ and let $K_m$ be the field obtained by adjoining $T^{1/q^m}$ to $K$. Then $K_m$ is a complete valued field with respect to the extended absolute value of $K$. Indeed, it is easy to check that $K_m = \mathbf{F}_q((T^{1/q^m}))$.

**Proposition 2** *Let $X$ be a compact subset of $K$, or more generally any compact Hausdorff totally disconnected topological space and let $\{f_n\}_{n\geq 0}$ is an orthonormal basis of $C(X, K)$. Then*
*(1) $\{f_n^{q^m}\}_{n\geq 0}$ is an orthonormal basis for $C(X, K)$ for all positive integer $m$.*
*(2) $\{f_n^{1/q^m}\}_{n\geq 0}$ is an orthonormal basis for $C(X, K_m)$ for all positive integer $m$.*

**Proof.** For Part (1), assuming $\{f_n\}_{n\geq 0}$ is an orthonormal basis of $C(X, K)$, we see that each $f_n$ maps $X$ into $O$, as part of the definition of an orthonormal basis. For any function $f : X \to O$, $f^q(x) \equiv f(x) \bmod(T)$, so the difference $f^q(x) - f(x)$ takes its values in the maximal ideal $(T)$ of $O$. Therefore $\|f^q - f\| \leq |T| < 1$. This upper bound is independent of $f$. Apply it to the functions $f_n(x)$ and invoke Lemma 1. Of course $q$ can be replaced with $q^m$, where $m$ is any positive integer.

Part (2) follows from the same argument above with $q$ replaced by $1/q$ as long as we know that $\{f_n\}_{n\geq 0}$ is also an orthonormal basis of $C(X, K_m)$. To see this, raising to $1/q^m$ power is an isomorphism from $K$ to the field $K_m$, which preserves the property of having norm $< 1$, equal to 1, and $> 1$. Then it is not difficult to check that this map carries an orthonormal basis of $C(O, K)$ to that of $C(O, K_m)$. □

**Corollary 5** *For each integer $m \geq 0$, both $\{G_j^{q^m}(x)\}_{j\geq 0}$ and $\{\mathsf{D}_j^{q^m}(x)\}_{j\geq 0}$ are orthonormal bases for $C(O, K)$.*

**Proof.** It follows from Theorem 7 and Proposition 2. □

**Corollary 6** *For any integer $m \geq 0$, both $\{E_n^{q^m}(x)\}_{n\geq 0}$ and $\{\mathcal{D}_n^{q^m}(x)\}_{n\geq 0}$ are orthonormal bases for $LC(O, K)$.*

**Proof.** It follows by applying the proof of Proposition 2 to $LC(O, K)$ instead of $C(O, K)$, together with Corollary 3. For a direct proof of $\{E_n^{q^m}(x)\}_{n\geq 0}$, one can deduce the equality

$$E_n^q(x) = [n+1]E_{n+1}(x) + E_n(x)$$

by using the recursion formula (see p. 46 [Go2]): $e_{n+1}(x) = e_n(x)^q - F_n^{q-1}e_n(x)$ and properties of $F_n$. Thus on taking norm of the equation we see that $\|E_n^q(x) - E_n(x)\| \leq 1/q < 1$. Now the case of $m = 1$ follows at once from Lemma 1. Generally for $m > 1$, we use the non-archimedean metric inequality to show that the aforementioned inequality with $q$ replaced by $q^m$ is also true. See Remark 2 in section 5 for an alternate proof of $\{\mathcal{D}_n^{q^m}\}_{n\geq 0}$. □

## 4 Properties of two bases

In this section, we will compare some properties of the Carlitz polynomials and digit derivatives when compared to their constructions. Carlitz developed some properties of Carlitz polynomials such as the addition law(binomial theorem)in the following.

**Proposition 3** (1) $G_j(\alpha x) = \alpha^j G_j(x)$ for $\alpha \in \mathbf{F}_q^*$.
(2) $G_j(x+u) = \sum_{e+f=j} \binom{j}{e} G_e(x) G_f(u)$.
(3) $G_j'(\alpha x) = \alpha^j G_j'(x)$ for $\alpha \in \mathbf{F}_q^*$.
(4) $G_j'(x+u) = \sum_{e+f=j} \binom{j}{e} G_e(x) G_f'(u)$.

**Proof.** See [C2] and [Go1] for a detailed proof. □

Similarly $\mathsf{D}_j$ also satisfy the same properties as the Carlitz polynomials do.

**Proposition 4** (1) $\mathsf{D}_j(\alpha x) = \alpha^j \mathsf{D}_j(x)$ for $\alpha \in \mathbf{F}_q^*$.
(2) $\mathsf{D}_j(x+u) = \sum_{e+f=j} \binom{j}{e} \mathsf{D}_e(x) \mathsf{D}_f(u)$.
(3) $\mathsf{D}_j'(\alpha x) = \alpha^j \mathsf{D}_j'(x)$ for $\alpha \in \mathbf{F}_q^*$.
(4) $\mathsf{D}_j'(x+u) = \sum_{e+f=j} \binom{j}{e} \mathsf{D}_e(x) \mathsf{D}_f'(u)$.

**Proof.** The proof is exactly the same as those given in [C2] and [Go1]. □

As $\binom{q^m-1}{\alpha} = (-1)^\alpha$ in $\mathbf{F}_q$, we get the following corollaries to Propositions 3 and 4.

**Corollary 7** (1) $G_{q^m-1}(x+u) = \sum_{e+f=q^m-1} (-1)^e G_e(x) G_f(u)$.
(2) $G_{q^m-1}(x-u) = \sum_{e+f=q^m-1} G_e(x) G_f(u)$.
(3) $G_{q^m-1}'(x+u) = \sum_{e+f=q^m-1} (-1)^e G_e(x) G_f'(u)$.
(4) $G_{q^m-1}'(x-u) = \sum_{e+f=q^m-1} G_e(x) G_f'(u)$.

**Corollary 8** (1) $\mathsf{D}_{q^m-1}(x+u) = \sum_{e+f=q^m-1} (-1)^e \mathsf{D}_e(x) \mathsf{D}_f(u)$.
(2) $\mathsf{D}_{q^m-1}(x-u) = \sum_{e+f=q^m-1} \mathsf{D}_e(x) \mathsf{D}_f(u)$.
(3) $\mathsf{D}_{q^m-1}'(x+u) = \sum_{e+f=q^m-1} (-1)^e \mathsf{D}_e(x) \mathsf{D}_f'(u)$.
(4) $\mathsf{D}_{q^m-1}'(x-u) = \sum_{e+f=q^m-1} \mathsf{D}_e(x) \mathsf{D}_f'(u)$.

Via an auxiliary function, Carlitz [C2] used polynomial interpolations to prove the orthogonality property of $G_j(x)$ and $G_j'(x)$. Recently Yang [Y] established the same result in an elementary way.

**Theorem 8** (1) *For $l < q^n$, $k$ an arbitrary integer $\geq 0$,*

$$\sum_{\substack{m \in \mathbf{F}_q[T] \\ deg(m) < n}} G_k(m) G_l'(m) = \begin{cases} 0 & if\, k+l \neq q^n - 1 \\ (-1)^n & if\, k+l = q^n - 1 \end{cases}$$

(2) *For $l < q^n$, $k < q^n$,*

$$\sum_{\substack{m \ monic \\ deg(m) = n}} G_k(m) G_l'(m) = \begin{cases} 0 & if\, k+l \neq q^n - 1 \\ (-1)^n & if\, k+l = q^n - 1 \end{cases}$$

Now we also have the orthogonality property of the digit derivatives parallel to Theorem 8.

**Theorem 9** (1) *For $l < q^n$, $k$ an arbitrary integer $\geq 0$,*

$$\sum_{\substack{m \in \mathbf{F}_q[T] \\ deg(m) < n}} \mathsf{D}_k(m) \mathsf{D}_l'(m) = \begin{cases} 0 & if\, k+l \neq q^n - 1 \\ (-1)^n & if\, k+l = q^n - 1 \end{cases}$$

(2) For $l < q^n, k < q^n$,

$$\sum_{\substack{m \text{ monic} \\ \deg(m) = n}} \mathsf{D}_k(m)\mathsf{D}'_l(m) = \begin{cases} 0 & if\, k+l \neq q^n - 1 \\ (-1)^n & if\, k+l = q^n - 1 \end{cases}$$

**Proof.** The proof of Theorem 1.9 [Y] is carried over word for word with $G_k(x), H_i^s := E_i(T^{i+s})$ replaced by $\mathsf{D}_k(x), \mathcal{D}_i(T^{i+s})$ respectively. We mention here that Yang's result works with the digit principle once we have a basis $\{F_i\}$ of $LC(O, K)$ provided that $F_i(T^e) = 0$ for $e < i$ and equals 1 for $e = i$. (Both $E_i$ and $\mathcal{D}_i$ satisfy this assumption.) □

Using the properties of $G_j$ in Proposition 3 Wagner characterized continuous $\mathbf{F}_q$-linear operators on $O$ in terms of expansion coefficients for a given expansion of $f \in C(O, K)$ with respect to $G_j$.

**Theorem 10** *Write $f \in C(O, K)$ as $f(x) = \sum_{j=0}^{\infty} A_j G_j(x)$. Then $f \in LC(O, K)$ if and only if $A_j = 0$ for $j \neq q^n$, where $n \geq 0$.*

**Proof.** It is obvious from Theorem 1 (Corollary 3) but for an alternate proof we refer to [W1] or [Y]. □

Similarly we can employ the identities in Proposition 4 to characterize continuous $\mathbf{F}_q$-linear operators on $O$ in terms of coefficients of the expansion of functions with respect to the digit derivatives in the following.

**Theorem 11** *Write $f \in C(O, K)$ as $f(x) = \sum_{j=0}^{\infty} B_j \mathsf{D}_j(x)$. Then $f \in LC(O, K)$ if and only if $B_j = 0$ for $j \neq q^n$, where $n \geq 0$.*

**Proof.** We see the result is obvious from Theorem 2 but without using it we here go along with the lines of the arguments in [W1].
Suppose $B_j = 0$ for $j \neq q^n$ where $n \geq 0$. Then the expansion of $f$ becomes

$$f(x) = \sum_{n=0}^{\infty} B_{q^n} \mathcal{D}_n(x).$$

Since the partial sums of the preceding series are $\mathbf{F}_q$-linear operators, we get the desired result.
Conversely, we use two identities (1) and (2) in Proposition 4. Let $\alpha$ be a primitive root of unity in $\mathbf{F}_q$. Then from (1) in Proposition 4 we get

$$\sum_{j=0}^{\infty} \alpha B_j \mathsf{D}_j(x) = \sum_{j=0}^{\infty} \alpha^j B_j \mathsf{D}_j(x),$$

so $B_j = 0$ unless $j \equiv 1 \bmod (q-1)$. Using the linearity of $f$ and (2), we get

$$\sum_{e=0}^{\infty} B_e \mathsf{D}_e(x) + \sum_{e=0}^{\infty} B_e \mathsf{D}_e(u) = \sum_{e=0}^{\infty} \{\sum_{j=e}^{\infty} \binom{j}{e} \mathsf{D}_{j-e}(u)\} \mathsf{D}_e(x).$$

Equating constant terms, i.e., $\mathsf{D}_0$, of both sides of the preceding equation, we see that $B_0 = 0$. Equating coefficients of $\mathsf{D}_e(x)$ for $e > 0$ and subtracting $B_e$, we get

$$\sum_{j=e+1}^{\infty} B_j \binom{j}{e} \mathsf{D}_{j-e}(u) = 0.$$

Hence, for all $e, j$ with $1 \leq e < j$,

$$\binom{j}{e} B_j = 0.$$

It yields that $B_j = 0$ unless $j = p^t$, where $p$ is the characteristic of $\mathbf{F}_q$. Since $p^t \equiv 1 \bmod (q-1)$, we must have $p^t = q^n$, where $n \geq 0$. □

## 5 Expansion of functions in $LC(O, K)$

From Corollary 6 we know for any non-negative integer $m$ that the set of $q^m$ powers of hyper-differential operators is an orthonormal basis of $LC(O, K)$. In this section, we will focus on obtaining the closed formula for expansion coefficients of any $f \in LC(O, K)$ with respect to $\{\mathcal{D}_n^{q^m}\}_{n \geq 0}$. Hence we generalize Theorem 2.

**Theorem 12** *Every $f \in LC(O, K)$ can be expanded uniquely as*

$$f(x) = \sum_{n=0}^{\infty} \beta_n^{(m)} \mathcal{D}_n^{q^m}(x),$$

*where the coefficients $\beta_n^{(m)}, n \geq 0$ are given by the formula*

$$\beta_n^{(m)} = (\Delta - [m]I)^n f(1) = \sum_{i=0}^{n} \sum_{j=i}^{n} (-1)^{n-i} \binom{n}{j} [m]^{n-j} f(T^i) \mathcal{D}_i(T^j).$$

**Proof.** We will use the properties of iterative derivations to show the coefficients of expansion of any $f \in LC(O, K)$ are uniquely recovered by the formula related to the Carlitz difference operator $\Delta$. Let $f = \sum_{i=0}^{\infty} \beta_i^{(m)} \mathcal{D}_i^{q^m}$ be a uniformly convergent series with respect to $\{\mathcal{D}_n^{q^m}\}_{n \geq 0}$. Then using the product rule (2) of $\mathcal{D}_n$, we easily see that for each $i \geq 1$,

$$\Delta(\mathcal{D}_i^{q^m}(x)) = [m]\mathcal{D}_i^{q^m}(x) + \mathcal{D}_{i-1}^{q^m}(x).$$

Hence we have $(\Delta - [m]I)\mathcal{D}_i^{q^m}(x) = \mathcal{D}_{i-1}^{q^m}(x)$ where $I$ is the identity operator on $LC(O, K)$. This implies that by iterating $(\Delta - [m]I)$ we get $(\Delta - [m]I)^n \mathcal{D}_i^{q^m}(x) = \mathcal{D}_{i-n}^{q^m}$. Then substitution 1 for $x$ in the expansion of $f$ gives us that for each $n \geq 0$,

$$\beta_n^{(m)} = ((\Delta - [m]I)^n f)(1).$$

Now one can explicitly calculate $\beta_n^{(m)} = ((\Delta - [m]I)^n f)(1)$ in the closed form as in the case of $m = 0$ (Theorem 2). We will here deal with these coefficients more generally. Let $g(y) = f(xy)$ for any continuous $\mathbf{F}_q$-linear function $f$ on $O$. Write $g(y) \in LC(O, K)$ as $g(y) = \sum_{n=0}^{\infty} \beta_n^{(m)}(x) \mathcal{D}_n^{q^m}(y)$, then the coefficients are given by

$$\beta_n^{(m)}(x) = ((\Delta - [m]I)^n g)(1) = ((\Delta - [m]I)^n f)(x).$$

By the binomial expansion, we get

$$\beta_n^{(m)}(x) = \sum_{j=0}^{n} \binom{n}{j} (-1)^{n-j} [m]^{n-j} \Delta^j f(x).$$

From the expansion of $f$ with respect to $\{\mathcal{D}_n\}$, it is known in [J] or [S] that

$$\Delta^j f(x) = \sum_{i=0}^{j} (-1)^{j-i} f(T^i x) \mathcal{D}_i(T^j).$$

Thus one can get

$$\beta_n^{(m)}(x) = \sum_{j=0}^{n} \sum_{i=0}^{j} (-1)^{n-i} \binom{n}{j} [m]^{n-j} f(T^i x) \mathcal{D}_i(T^j).$$

Hence, as $\mathcal{D}_i(T^j) = 0$ for $i > j$,

$$\beta_n^{(m)}(x) = \sum_{i=0}^{n} \sum_{j=i}^{n} (-1)^{n-i} \binom{n}{j} [m]^{n-j} f(T^i x) \mathcal{D}_i(T^j).$$

Finally substituting 1 for $x$ in the equation above we have the coefficients

$$\beta_n^{(m)} := \beta_n^{(m)}(1) = \sum_{i=0}^{n}\sum_{j=i}^{n}(-1)^{n-i}\binom{n}{j}[m]^{n-j}f(T^i)\mathcal{D}_i(T^j).$$

Thus we complete the proof. □

The discussion above in the proof concludes the following.

**Corollary 9** *For any $f \in LC(O, K)$ and any integer $m \geq 0$, we have*

$$((\Delta - [m]I)^n f)(x) = \sum_{i=0}^{n}\sum_{j=i}^{n}(-1)^{n-i}\binom{n}{j}[m]^{n-j}f(T^i x)\mathcal{D}_i(T^j).$$

**Corollary 10** *Let $f(x) = \sum_{n=0}^{\infty}\beta_n^{(m)}\mathcal{D}_n^{q^m}(x)$ be a continuous $\mathbf{F}_q$-linear function from $O$ to $K$. Then $f \in LC(O, O)$ if and only if $\{\beta_n^{(m)}\}_{n\geq 0} \subset O$.*

It is interesting to find the explicit expansion of $\mathcal{D}_n^{q^m}(x)$ in terms of $\mathcal{D}_n(x)$ and vice versa.

**Proposition 5** (1) *For each $n \geq 0$ and $x \in O, \mathcal{D}_n^{q^m}(x) = \sum_{i=0}^{\infty}[m]^i\binom{i+n}{n}\mathcal{D}_{i+n}(x)$.*
(2) *For each $n \geq 0$ and $x \in O, \mathcal{D}_n(x) = \sum_{i=0}^{\infty}(-[m])^i\binom{i+n}{n}\mathcal{D}_{i+n}^{q^m}(x)$.*

**Proof.** For Part (1) use a Voloch's identity [V]:

$$x^{q^m} = \sum_{i=0}^{\infty}[m]^i\mathcal{D}_i(x).$$

Substituting $\mathcal{D}_n(x)$ for $x$ into the equation, we use the composite rule (3) of $\mathcal{D}_n$ to obtain

$$\mathcal{D}_n(x)^{q^m} = \sum_{i=0}^{\infty}[m]^i\mathcal{D}_i(\mathcal{D}_n(x))$$
$$= \sum_{i=0}^{\infty}[m]^i\binom{i+n}{n}\mathcal{D}_{i+n}(x).$$

Alternatively, note that both sides of the equation (1) are $\mathbf{F}_q$-linear, so it suffices to check the formula for $x = T^s, s \geq 0$ and in this case it is straightforward. Now we use the preceding argument to see Part (2) holds. Substituting $T^s$ for $x$, we see the right-hand side of the equation (2) equals

$$\sum_{i=0}^{s-n}(-[m])^i\binom{i+n}{n}\binom{s}{i+n}T^{(s-n-i)q^m} = \binom{s}{n}\sum_{i=0}^{s-n}(-[m])^i\binom{s-n}{i}T^{(s-n-i)q^m}$$
$$= \binom{s}{n}(-[m]+T^{q^m})^{s-n},$$

which is equal to the left-hand side of the equation (2). Another way of showing Part (2) is to use a crucial identity (Proposition 5.2.16 in[J]) for continuous $\mathbf{F}_q$-linear functions on $O$:

$$x^{1/q^m} = \sum_{i=0}^{\infty}(T^{1/q^m} - T)^i\mathcal{D}_i(x).$$

Substitute $\mathcal{D}_n(x)$ for $x$ into the above equation, then the composite rule of $\mathcal{D}_n$ yields

$$\mathcal{D}_n(x)^{1/q^m} = \sum_{i=0}^{\infty}(T^{1/q^m} - T)^i\mathcal{D}_i(\mathcal{D}_n(x))$$
$$= \sum_{i=0}^{\infty}(T^{1/q^m} - T)^i\binom{i+n}{n}\mathcal{D}_{i+n}(x).$$

By taking $q^m$th power of both sides, we get

$$\begin{aligned}\mathcal{D}_n(x) &= \sum_{i=0}^{\infty}(T-T^{q^m})^i\binom{i+n}{n}\mathcal{D}_{i+n}^{q^m}(x) \\ &= \sum_{i=0}^{\infty}(-[m])^i\binom{i+n}{n}\mathcal{D}_{i+n}^{q^m}(x).\end{aligned}$$

$\square$

**Remarks**
1. Imitating the arguments of Wagner([W2]), together with the equation in Corollary 9, one can also prove the same result in Theorem 12 by showing that for a continuous $\mathbf{F}_q$-linear function $f$, the series $\sum_{n=0}^{\infty}((\Delta-[m]I)^n f)(1)\mathcal{D}_n^{q^m}(x)$ converges uniformly to $f$ on the valuation ring of $K$.
2. Note that Proposition 5(2) is obtained without using the fact that $\{\mathcal{D}_n^{q^m}\}_{n\geq 0}$ is an orthonormal basis of $LC(O,K)$, so one also give an alternate proof for $\mathcal{D}_n^{q^m}$ in Corollary 6 by showing that for each $n\geq 0$,

$$\|\mathcal{D}_n^{q^m}-\mathcal{D}_n\| = \max_{i\geq 1}\{|-[m]^i\binom{i+n}{n}|\} \leq 1/q < 1.$$